\documentclass[12pt,amstex,amscd,graphicx,euscript]{amsart}

\begin{document}

\long\def\state#1#2{
\medskip\par\noindent
{\bf #1} 
{\it #2}
\par\medskip
}

\newtheorem{theorem}{Theorem}[section]
\newtheorem{prop}[theorem]{Proposition}
\newtheorem{lemma}[theorem]{Lemma}
\newtheorem{cor}[theorem]{Corollary}
\newtheorem{defn}[theorem]{Definition}
\newtheorem{conj}[theorem]{Conjecture}
\newtheorem{claim}[theorem]{Claim}
\newtheorem{remark}[theorem]{Remark}
\newtheorem{qn}[theorem]{Question}

\newcommand{\gothic}{\mathfrak}

\newcommand{\bdy}{\partial}

\newcommand{\bA}{{\mathbb A}}
\newcommand{\bB}{{\mathbb B}}
\newcommand{\bC}{{\mathbb C}}
\newcommand{\bD}{{\mathbf D}}
\newcommand{\bE}{{\mathbb E}}
\newcommand{\bF}{{\mathbf F}}
\newcommand{\bG}{{\mathbf G}}
\newcommand{\bH}{{\mathbf H}}
\newcommand{\bI}{{\mathbf I}}
\newcommand{\bJ}{{\mathbf J}}
\newcommand{\bK}{{\mathbf K}}
\newcommand{\bL}{{\mathbb L}}
\newcommand{\bM}{{\mathbb M}}
\newcommand{\bN}{{\mathbb N}}
\newcommand{\bO}{{\mathbb O}}
\newcommand{\bP}{{\mathbb P}}
\newcommand{\bQ}{{\mathbb Q}}
\newcommand{\bR}{{\mathbb R}}
\newcommand{\bS}{{\mathbf S}}
\newcommand{\bT}{{\mathbf T}}
\newcommand{\bU}{{\mathbf U}}
\newcommand{\bV}{{\mathbf V}}
\newcommand{\bW}{{\mathbf W}}
\newcommand{\bX}{{\mathbf X}}
\newcommand{\bY}{{\mathbf Y}}
\newcommand{\bZ}{{\mathbb Z}}

\newcommand{\calA}{{\mathcal A}}
\newcommand{\calB}{{\mathcal B}}
\newcommand{\calC}{{\mathcal C}}
\newcommand{\calD}{{\mathcal D}}
\newcommand{\calE}{{\mathcal E}}
\newcommand{\calF}{{\mathcal F}}
\newcommand{\calG}{{\mathcal G}}
\newcommand{\calH}{{\mathcal H}}
\newcommand{\calI}{{\mathcal I}}
\newcommand{\calJ}{{\mathcal J}}
\newcommand{\calK}{{\mathcal K}}
\newcommand{\calL}{{\mathcal L}}
\newcommand{\calM}{{\mathcal M}}
\newcommand{\calN}{{\mathcal N}}
\newcommand{\calO}{{\mathcal O}}
\newcommand{\calP}{{\mathcal P}}
\newcommand{\calQ}{{\mathcal Q}}
\newcommand{\calR}{{\mathcal R}}
\newcommand{\calS}{{\mathcal S}}
\newcommand{\calT}{{\mathcal T}}
\newcommand{\calU}{{\mathcal U}}
\newcommand{\calV}{{\mathcal V}}
\newcommand{\calW}{{\mathcal W}}
\newcommand{\calX}{{\mathcal X}}
\newcommand{\calY}{{\mathcal Y}}
\newcommand{\calZ}{{\mathcal Z}}

\newcommand\CH{{\calC\calH}}
\newcommand\MF{{\calM\calF}}
\newcommand\PMF{{\calP\kern-2pt\calM\calF}}
\newcommand\ML{{\calM\calL}}
\newcommand\PML{{\calP\kern-2pt\calM\calL}}
\newcommand\GL{{\calG\calL}}

\newcommand{\gA}{{\mathfrak A}}
\newcommand{\gB}{{\mathfrak B}}
\newcommand{\gC}{{\mathfrak C}}
\newcommand{\gD}{{\mathfrak D}}
\newcommand{\gE}{{\mathfrak E}}
\newcommand{\gF}{{\mathfrak F}}
\newcommand{\gG}{{\mathfrak G}}
\newcommand{\gH}{{\mathfrak H}}
\newcommand{\gI}{{\mathfrak I}}
\newcommand{\gJ}{{\mathfrak J}}
\newcommand{\gK}{{\mathfrak K}}
\newcommand{\gL}{{\mathfrak L}}
\newcommand{\gM}{{\mathfrak M}}
\newcommand{\gN}{{\mathfrak N}}
\newcommand{\gO}{{\mathfrak O}}
\newcommand{\gP}{{\mathfrak P}}
\newcommand{\gQ}{{\mathfrak Q}}
\newcommand{\gR}{{\mathfrak R}}
\newcommand{\gS}{{\mathfrak S}}
\newcommand{\gT}{{\mathfrak T}}
\newcommand{\gU}{{\mathfrak U}}
\newcommand{\gV}{{\mathfrak V}}
\newcommand{\gW}{{\mathfrak W}}
\newcommand{\gX}{{\mathfrak X}}
\newcommand{\gY}{{\mathfrak Y}}
\newcommand{\gZ}{{\mathfrak Z}}

\newcommand{\eA}{{\EuScript A}}
\newcommand{\eB}{{\EuScript B}}
\newcommand{\eC}{{\EuScript C}}
\newcommand{\eD}{{\EuScript D}}
\newcommand{\eE}{{\EuScript E}}
\newcommand{\eF}{{\EuScript F}}
\newcommand{\eG}{{\EuScript G}}
\newcommand{\eH}{{\EuScript H}}
\newcommand{\eI}{{\EuScript I}}
\newcommand{\eJ}{{\EuScript J}}
\newcommand{\eK}{{\EuScript K}}
\newcommand{\eL}{{\EuScript L}}
\newcommand{\eM}{{\EuScript M}}
\newcommand{\eN}{{\EuScript N}}
\newcommand{\eO}{{\EuScript O}}
\newcommand{\eP}{{\EuScript P}}
\newcommand{\eQ}{{\EuScript Q}}
\newcommand{\eR}{{\EuScript R}}
\newcommand{\eS}{{\EuScript S}}
\newcommand{\eT}{{\EuScript T}}
\newcommand{\eU}{{\EuScript U}}
\newcommand{\eV}{{\EuScript V}}
\newcommand{\eW}{{\EuScript W}}
\newcommand{\eX}{{\EuScript X}}
\newcommand{\eY}{{\EuScript Y}}
\newcommand{\eZ}{{\EuScript Z}}

\newcommand\half{{\textstyle{\frac12}}}
\newcommand\Half{{\frac12}}
\newcommand\ep{\epsilon}
\newcommand\hhat{\widehat}
\newcommand\til{\widetilde}
\newcommand\gesim{\succ}
\newcommand\lesim{\prec}
\newcommand\simle{\lesim}
\newcommand\simge{\gesim}
\newcommand{\simmult}{\asymp}
\newcommand{\simadd}{\mathrel{\overset{\text{\tiny $+$}}{\sim}}}
\newcommand{\sm}{\setminus}
\newcommand{\pair}[1]{\langle #1\rangle}
\newcommand{\tprec}{\prec_t}
\newcommand{\fprec}{\prec_f}
\newcommand{\bprec}{\prec_b}

\newcommand{\sprec}{\prec_s}
\newcommand{\cpreceq}{\preceq_c}
\newcommand{\cprec}{\prec_c}
\newcommand{\topprec}{\prec_{\rm top}}
\newcommand{\Topprec}{\prec_{\rm TOP}}
\newcommand{\fsub}{\mathrel{\scriptstyle\searrow}}
\newcommand{\bsub}{\mathrel{\scriptstyle\swarrow}}
\newcommand{\fsubd}{\mathrel{{\scriptstyle\searrow}\kern-1ex^d\kern0.5ex}}
\newcommand{\bsubd}{\mathrel{{\scriptstyle\swarrow}\kern-1.6ex^d\kern0.8ex}}
\newcommand{\fsubeq}{\mathrel{\raise-.7ex\hbox{$\overset{\searrow}{=}$}}}
\newcommand{\bsubeq}{\mathrel{\raise-.7ex\hbox{$\overset{\swarrow}{=}$}}}
\newcommand{\bbar}{\overline}

\newcommand\AM{{\mathcal{M}^{\omega}}}
\newcommand\AP{{\mathcal{P}^{\omega}}}
\newcommand\AC{{\mathcal{C}^{\omega}}}
\newcommand\AQ{{\mathcal{Q}^{\omega}}}
\newcommand\ACx{{\mathcal{C}^{\omega}}_\xi}

\newcommand\homeo{\cong}
\newcommand\indhat{{\hhat{\operatorname{ind}}}}
\newcommand\ind{{\operatorname{ind}}}
\newcommand\Ind{{\operatorname{Ind}}}
\newcommand\Dim{{\operatorname{dim}}}
\newcommand\Indhat{{\hhat{\operatorname{Ind}}}}
\newcommand\dimhat{{\hhat{\operatorname{dim}}}}
\def\MCG{\mathcal {MCG}}
\newcommand\tint[2]{#1 \pitchfork #2 \ne \emptyset}
\newcommand\notint[2]{#1 \pitchfork #2 = \emptyset}
\def\ulim{\lim_\omega}
\def\ulimf{\lim_\omega \frac{1}{s_i}}
\def\cone{{\rm{Cone}}_{\omega}}
\def\dist{{\rm{dist}}}
\def\co{\colon}

\newcommand\oseq[1]{\mathbf{#1}}
\newcommand\seq[1]{\mbox{\boldmath$#1$}}

\title{Mapping Class Groups and Interpolating Complexes: Rank }

\author{Mahan Mj.}
\address{RKM Vivekananda University, Belur Math, WB-711 202, India}

\date{}

\begin{abstract} 
A family of interpolating graphs $\calC (S, \xi )$ of complexity $\xi$ is constructed for a surface $S$ and $-2 \leq \xi \leq \xi (S)$. For $\xi = -2, -1, \xi (S) -1$ these specialise to graphs quasi-isometric to the marking graph, the pants graph and the curve graph respectively. We generalise  Theorems of Brock-Farb and Behrstock-Minsky to show that the rank of $\calC (S, \xi )$ is $r_\xi$, the largest number of disjoint copies of subsurfaces of complexity greater than $\xi $ that may be embedded in $S$. The interpolating graphs 
$\calC (S, \xi )$  interpolate between the pants graph and the curve graph. 

\begin{center}
AMS subject classification =   20F67(Primary), 22E40  
\end{center}

\end{abstract}

\maketitle

\tableofcontents

\section{Introduction}

\subsection{Motivation and Statement of Results}
Starting with Masur-Minsky's result that the curve complex is hyperbolic \cite{masur-minsky}, the coarse geometry of mapping class groups has attracted much attention. A motivating scholium is the following. \\

\medskip

\noindent {\it The Mapping Class Group behaves like a non-uniform rank one lattice away from the peripheral subgroups and like a higher rank non-uniform lattice at the peripheral subgroups.}

\medskip

Hyperbolicity of the curve complex is an instance of rank one behaviour. Intersection patterns of peripheral subgroups and resultant structures similar to the Tits complex illustrate higher rank behaviour.
In this paper, we investigate further this interbreeding of rank one and higher rank behaviour. 

There are two pieces of motivation behind this paper:\\
{\bf Motivation 1:} Higher rank lattices admit a whole family of compactifications, for instance the Borel-Serre, Reduced Borel-Serre and toroidal compactifications.  These are built from configurations of parabolic subgroups. (See Borel-Ji \cite{borel-ji} for instance.) Of particular relevance to this paper is the fact that the Furstenberg (or maximal) boundary is obtained as a quotient space of the Tits boundary by identifying certain Weyl chambers at infinity to points. 
{\bf A  coarse geometric analog of such a topological quotienting operation is  "coning" } (see below). This intuitive idea will play an important role in the construction of interpolating graphs.

 Moduli spaces too admit such compactifications, of which the Deligne-Mumford compactification is probably the most well-known. If we look at the universal cover of the compactified moduli space, we find an intersection pattern of boundary strata. This is encoded in the Curve Complex $\calC (S)$ of a surface $S$ of (discovered by Harvey \cite{harvey-cc} from the above-mentioned analogy with non-uniform lattices of higher rank). Another such simplicial complex (originally dicovered by Hatcher and Thurston \cite{hatcher-thurston} ) is the pants complex $\calP (S)$. Recently, Brock \cite{brock-jams} has shown that the pants complex is quasi-isometric to Teichmuller space equipped with the Weil-Petersson metric. The first aim of
this paper is to describe a collection of simplicial complexes interpolating (in a sense to be made precise) between the curve-complex $\calC (S)$ and the pants complex $\calP (S)$.

\medskip

\noindent {\bf Motivation 2:} Masur and Minsky develop in \cite{masur-minsky2} a detailed combinatorial structure of {\em hierarchies} to get a handle on quasigeodesics in $MCG(S)$. We develop in this paper a related hierarchy of spaces where the bottom level is given by the curve complex and the top level by the marking complex (quasi-isometric to the the mapping class group). Thus in a sense again, we describe a collection of simplicial complexes interpolating  between the curve-complex $\calC (S)$ and the marking complex $\calM (S)$.

\medskip

\noindent {\bf Quasi-isometric Models} We first describe quasi-isometric models for these interpolating complexes as they are easy to define. Let $\Gamma (S)$ be a Cayley graph of the mapping class group of a surface $S$ of genus $g$ with $n$ punctures. Let $\xi (S ) = 3g - 3 + n = \xi_0$ denote the complexity of the surface $S$. Note that $\Gamma (S)$ depends on the choice of generators and is therefore well-defined upto quasi-isometry. 

We now describe a variety of graphs associated with $\Gamma (S)$. Fix a $\xi$ with $-2 \leq \xi \leq \xi_0 - 1 $.

Consider all essential (i.e. $\pi_1$-injective) subsurfaces $S$ of complexity less than or equal to $\xi$. These fall into finitely many orbits under the action of the mapping class group, $S_1, \cdots S_k$, say. Then we may assume (changing the generating set if necessary) that the Cayley graph $\Gamma_{MCG(S_i)} \subset \Gamma (S)$. We may further assume for convenience that the $S_i$'s are maximal (i.e. no $S_i$ is a proper subset of another $S_j$.) We define $\Gamma (S, \xi )$ to be the graph obtained from $\Gamma (S)$ by coning (a la Farb \cite{farb-relhyp}) cosets of sub-mapping class groups $MCG(S_i)$, i.e. by introducing a vertex $v_{gH}$ for each coset $gH$ of the mapping class group $H=MCG(S_i)$ for each of the above $S_i$'s and joining it by an edge of length $\half$ to every element of $gH \subset \Gamma (S)$. The main theorem of this paper  determines the rank of $\Gamma (S, \xi )$.  In the next subsection, we shall give a more intrinsic (geometric) model $\calC (S, \xi )$ quasi-isometric to $\Gamma (S, \xi )$ with the additional restriction that it is defined {\em naturally}. These shall be termed the {\em interpolating graphs} or {\em complexity $\xi$ graphs} for the surface $S$. Note that $\xi = -1$ gives a graph qi to the pants graph and $\xi = -2$ cones off the trivial sub-mapping class groups corresponding to mapping class groups of disks, which are trivial, yielding therefore the mapping class group (or equivalently the marking graph). 

Let $r_\xi$ denote the maximum number of disjoint subsurfaces of complexity $(\xi + 1)$ that can be embedded in $S$. Then the main theorem of this paper states:

\medskip

\noindent {\bf Theorem \ref{rank} :}
{\it The rank of the interpolating graph $\calC (S, \xi )$ of complexity $\xi$ or its quasi-isometric model  $\Gamma (S, \xi )$ is $r_\xi   (S ) $. }

\medskip

In the final subsection of this paper, we shall draw a conjectural picture of the interconnections between hierarchies, rank and interpolating graphs. One interesting fallout is the following.
 
\medskip

\noindent {\bf Conjecture: Rank one implies Hyperbolic}\\
If for $\xi > 0$,  $r_\xi (S) = 1$  
 is $\calC (S, \xi )$ a hyperbolic metric space?

\subsection{Complexes Associated to the Mapping Class Group}

In this subsection, we describe the motivating examples: the curve-graph, the pants graph and the marking graph. We then proceed to give the promised description of interpolating (or complexity $\xi$) graphs.

\medskip

\noindent {\bf Curve graph}

\noindent {\bf Case 1: $\xi( S  ) \geq 2$} \\ 
The 
{\it curve graph} of $S$, denoted $\calC ( S  )$
is a graph with \\
1) vertices corresponding 
to nontrivial homotopy classes of non-peripheral, simple closed curves on $S$ \\
2) 
edges corresponding 
to pairs of (homotopically distinct)
simple closed curves which can be realized 
disjointly on $S$. \\

\smallskip

\noindent {\bf Case 2: $\xi( S  ) = 1$} \\ 
Then either $g=n=1$ (
$S$ is a one-holed torus ) or $g =0, n= 3$ ($S$ is a 4-holed
sphere). The vertex set is as in Case 1 above. The edge set consists of pairs of curves which realize the minimal possible intersection on 
$S$ (1 for the one-holed torus, 2 for the 4-holed sphere). 

\medskip

\noindent {\bf Case 3: $\xi( S  ) = 0$} \\ 
In this case, $S$ is the 3-holed  sphere with empty curve-complex as the vertex set  is empty.

\medskip

\noindent {\bf Case 4: $\xi( S  ) = -1$} \\ 
In this case, $S$ is the 2-holed  sphere or
 an annulus. Fix a point in each boundary component.\\
1) Vertices are (homotopy classes of)
 paths connecting the give boundary points up to homotopy rel endpoints. \\
2)Edges are pairs of non-intersecting
paths.\\

\medskip

\noindent {\bf Pants Graph}

\noindent The {\it pants graph} $\calP   (S ) $ of 
$S$ is a graph with \\
\begin{enumerate}
\item Vertices consisting of pants decompositions
of $S$. \\
\item  Edges consist of pants decompositions 
that agree on all but one curve, and further, those 
curves differ by an edge in the curve complex of the 
complexity one subsurface (complementary to the rest of the curves) in which they lie.
\end{enumerate}    

\medskip

\noindent {\bf Marking Graph}

\noindent We consider pants decompositions $\mu$ of $S$ (We may identify $\mu$ with maximal simplices in $\calC ( S  )$.
The set of underlying curves shall be denoted {\it base}
$ (\mu )$.
The transversals $\tau_\mu$ of $\mu$ consist of one curve for
    each component of {\emph base} $(\mu)$, intersecting it transversely. A pair  ( {\emph base} $ (\mu), {\tau_\mu }$) of pants decompositions and transversals shall be referred to as a {\em marking}.

\smallskip

A marking  (base $(\mu), {\tau_\mu }$)
is {\it clean, complete}  if, \\
1) for each $\gamma \in {base}$ $(\mu)$, the
transversal curve $t_\gamma$ to $\gamma$ is disjoint from the rest of  
{\it base} $(\mu)$. \\
2) Each pair $(\gamma , t_\gamma )$ fills a non-annular surface $W$ satisfying 
$\xi (W)=1$ and for which $d_{\calC (W)}(\gamma, t)=1$.\\

\smallskip

The {\it marking graph} or {\it marking complex}, 
$\calM (S)$ is defined as follows. \\
Vertices correspond to clean complete markings. 
The edges of $\calM (S)$ are of two types (See Masur-Minsky \cite{masur-minsky2} and also Behrstock-Minsky \cite{behrstock-minsky}.): \\
 {\it Twist}: Replace a transversal curve by another 
    obtained by performing a Dehn twist along the associated base 
    curve. \\
{\it Flip}: Exchange
 the roles of a base curve and its
     transversal curve. Perform surgery if necessary to reinstate disjointness.  
    (Note that after switching base and transversal, the disjointness requirement on
    the transversals may be violated.  However, Masur-Minsky show in \cite{masur-minsky2} that one can surger the new transversal 
    to obtain one that does satisfy the disjointness requirement.
    They further show that only a finite (uniformly bounded) number of flip moves are possible.)

\medskip

\begin{lemma} {\bf Masur-Minsky} \cite{masur-minsky2} $\calM ( S  )$ is quasi-isometric to the mapping class 
    group of $S$.
\label{mcg=markingcx}
\end{lemma}

\begin{remark} \cite{masur-minsky2}
Masur and Minsky note that 
the pants graph is exactly what remains of the marking complex when 
annuli (and hence transverse curves) are ignored. 
\label{mcpc}
\end{remark}

\begin{center}

{\bf Interpolating Graphs}

\end{center}

\medskip

We are now in a position to define the (natural or) geometric models that are the main object of study in this paper. 

\begin{defn}
An {\em interpolating graph} or  {\em complexity $\xi$ graph} $\calC (S, \xi )$ consists of the following. \\
1) vertices are pants decompositions (or maximal simplices in the curve complex) \\
2) edges are of two types:\\
a) edges of the pants graph of $S$ \\
b) additional edges connecting
 pairs of pants decompositions agreeing on the complement of a (connected) subsurface of complexity less than or equal to $\xi$\\
\label{intgr}
\end{defn}

Note that edges of type (2b) above {\it include} edges of type (2a) for $ \xi \geq 1$. However, in order to include
the pants graph as a starting point, we have mentioned edges of type (2a) separately.

As in \cite{masur-minsky2}, the same definitions apply to essential (possibly disconnected)
subsurfaces of $S$. 
For a disconnected 
surface $\displaystyle W=\sqcup_{i=1}^{n} W_{i}$, 
$\calC (W, \xi )=\prod_{i=1}^{n}\calC (W_{i}, \xi )$.

\begin{remark}
{\rm We note that the  interpolating complexes $\calC (S, \xi )$ may be regarded as 
(quasi-isometric models) unifying the 3 types of complexes given above:\\
1) $\xi = \xi_0 -1$: (denoting $\xi (S)$ by $\xi_0$.)
Subsurfaces where moves are considered are (arbitrary) proper subsurfaces. This gives a model quasi-isometric to the curve-graph. \\
To see this, note that all proper pants subgraphs (i.e. pants graphs of all proper
essential subsurfaces) have diameter one in  $\calC (S, \xi )$ here. Recall (Remark \ref{mcpc})
that 
the pants graph is exactly what remains of the marking complex when 
annuli (and hence transverse curves) are ignored, or equivalently, when all moves on annulii are at distance one from
each other. Thus $\calC (S, \xi )$ is quasi-isometric to what one gets from the marking complex by first
forgetting annulii and then all subsurfaces of complexity $\xi$. This is equivalent to coning all proper 
mapping class subgroups (i.e. mapping class groups of all proper
essential subsurfaces) of the mapping class group of $S$, which is quasi-isometric to
the curve complex $CC(S)$ of $S$. \\
An explicit quasi-isometry from $CC(S)$ to $\calC (S, \xi )$ can be set up by sending any curve $\zeta \in 
CC(S)$ to some (any) pants decomposition containing $\zeta$. There is a bounded amount of ambiguity in this
as the set of pants decomposition containing $\zeta$ has diameter one in $\calC (S, \xi )$.\\
2) $\xi = -1$: This coincides exactly with the definition of the pants graph. There are no edges of type
(2b) of Definition \ref{intgr} (vacuously). We mention this case separately to underscore the point that
the pants graph is what one gets when annuli (and hence transverse curves) are ignored. \\
3) $\xi = -2$: Here, moves are restricted to surfaces of complexity $-2$, i.e. disks. To corectly interpret this case, the curve graph of the annulus has to be resurrected and re-instated over  the pants graph. This gives the marking graph. (Alternately, 
we could start off by defining our interpolating complex starting with the marking complex and then adjoin edges for moves occurring in subsurfaces of complexity $\leq \xi$. This takes care of the marking complex; but in all subsequent discussions, it would make the exposition more awkward - hence the definition here.)}
\label{ic}
\end{remark}

\begin{remark}

Given the definition above, it is easy to see that (for some, hence any, finite generating set) $\Gamma (S, \xi )$ is quasi-isometric to $\calC (S, \xi )$. 

In particular, (upto quasi-isometry) coning Dehn twists gives the pants complex and coning all proper sub-mapping class groups the curve complex (see \cite{masur-minsky} ). Not coning anything (or coning only the trivial sub-mapping class group for disks) gives the marking complex. 
\end{remark}

\begin{remark}
A useful heuristic is: \\
The complexity $\xi$ graph is what remains of the marking complex when surfaces of complexity $\leq \xi$ are ignored. 
\label{heuristic}
\end{remark}

\subsection{Projections, Hierarchies, Distance Formulae}

In this subsection, we summarise some of the foundational work of Masur-Minsky \cite{masur-minsky} \cite{masur-minsky2}, followed by more recent work of
Behrstock-Minsky \cite{behrstock-minsky}. An essential
tool for the next section
is Theorem \ref{dfic} giving a distance formula for intepolating complexes.

\begin{theorem}
{\bf Masur-Minsky \cite{masur-minsky} }
For any surface $S$, 
    the complex of curves $\calC   (S ) $ is an infinite diameter $\delta$-hyperbolic 
    space (as long as it is non-empty).
\label{cc-hyp}
\end{theorem}

\begin{defn} {\bf Masur-Minsky \cite{masur-minsky2}}
Given a subsurface $ W \subset  S $,  a 
\emph{subsurface 
projection}  is a map $\pi_W: \calC   (S )  \to 2^{\calC (W)}$ defined as follows. \\
Case 1:  $W$ is not an annulus. \\
Given any 
curve $\gamma \in \calC   (S ) $ intersecting $Y$ essentially, 
we define 
$\pi_W(\gamma)$ to be the collection of curves (vertices) in $\calC (W)$ 
obtained by surgering the essential arcs of $\gamma\cap W$ along 
$\partial W$ to obtain simple closed 
curves in $W$. \\
 Case 2:  $W$ is an annulus. \\
Here, the curve-graph is assumed to be that for  the compactified cover $\til W$ of $S$ corresponding to the subgroup $\pi_1(W)$. Note that $\til W$ can be identified with the annulus. 
If  $\gamma$ intersects $W$ transversely and
essentially, we 
lift $\gamma$ to an arc crossing the annulus $\til{W}$ and let this
be $\pi_W(\gamma)$. If $\gamma$ is a core curve of $W$ or fails
to intersect it, we let $\pi_W(\gamma) = \emptyset$.
\end{defn}

$d_{\calC (W)}(\mu,\nu)$ will be used as a short form for $d_{\calC (W)}(\pi_W(\mu),\pi_W(\nu))$.

Next,  for any $\mu\in\calC  (S, \xi ) $ and any non-annular $W\subseteq  S $ 
the above projection map induces  $\pi_W:\calC (S, \xi )\to 2^{\calC (W)}$. 
This map is simply the union over $\gamma\in$ {\bf base}$(\mu)$ of the 
usual projections $\pi_W(\gamma)$. As in the case of curve 
complex projections, we write $d_{\calC (W)}(\mu,\nu)$ for 
$d_{\calC (W)}(\pi_{ (W)}(\mu),\pi_{ (W)}(\nu))$. 
The distance in the interpolating graph of complexity $\xi$ shall be denoted by $d_\xi$.

\medskip

\noindent
 {\bf Hierarchies}

\noindent
As summarized in \cite{behrstock-minsky}, hierarchy paths are quasigeodesics in $\calM (S)$ with
constants depending only on the topological type of $S$ such that \\
\begin{enumerate}
\item any two points $\mu,\nu\in\calM (S)$ are connected by at least one
hierarchy path $\gamma$. \\
\item There is a monotonic map $v:\gamma\to\beta$, such that
$v(\gamma_n)$ is a vertex in {\bf base} $(\gamma_n)$ for every $\gamma_n$ in $\gamma$. \\
\item 
Subsurfaces of $S$ which ``separate'' $\mu$
from $\nu$ in a significant way must play a role in the hierarchy
paths from $\mu$ to $\nu$ in the following sense: \\
There exists a constant $M_{2}=M_{2}(S)$ such that, if 
$W$ is an essential
subsurface of $S$ and $d_{\calC (W)}(\mu,\nu) > M_{2}$, then for any
hierarchy path $\gamma$ connecting $\mu$ to $\nu$, there exists
a marking $\gamma_n$ in $\gamma$ with $[\partial W] \subset $
{\bf base} $(\gamma_n)$. Furthermore there exists a vertex $v$ in the geodesic
$\beta$ shadowed by $\gamma$ such that $W\subset S\setminus v$. 
This property follows directly from Lemma 6.2 of \cite{masur-minsky2}.
\end{enumerate}

\medskip

\noindent
 {\bf Distance Formulae}

\smallskip

Masur--Minsky prove the following distance formula for distances in the marking complex:\\
\begin{theorem}{\bf Masur--Minsky \cite{masur-minsky2}}
    If $\mu,\nu\in \calM (S)$,  then there exists a constant $K  (S ) $, 
    depending only on $S$, such that for each 
    $K>K  (S ) $ there exists $a\geq 1$ and $b\geq 0$ for which:
    $$
    d_{\calM (S)}(\mu,\nu) \approx_{a,b} \sum_{W\subseteq S} 
   Tf_K{d_{\calC (W)}(\pi_{Y}(\mu),\pi_{Y}(\nu))}
    $$
\label{dfmc}
\end{theorem}

Here the {\em threshold function}
$Tf_K{N}$ is defined to be $N$ if $N>K$ and $0$ else. 

In \cite{behrstock-minsky}, Behrstock and Minsky show that this yields the following formula for distances in the pants complex. They use Remark \ref{mcpc} to forget projections to annulii.

\begin{theorem}{\bf  Behrstock-Minsky \cite{behrstock-minsky} }
    If $\mu,\nu\in \calP (S)$,  then there exists a constant $K(S)$, 
    depending only on the topological type of $S$, such that for each 
    $K>K(S)$ there exists $a\geq 1$ and $b\geq 0$ for which:
    $$
    d_{\calP (S)}(\mu,\nu) \approx_{a,b} 
    \sum_{\rm{non-annular \,}Y\subseteq S} 
 Tf_K{d_{\calC (Y)}(\pi_{Y}(\mu),\pi_{Y}(\nu))}
    $$
\label{dfpc}
\end{theorem}

Exactly analogously, we state the generalisation of  Theorem \ref{dfpc} above to interpolating graphs 
$\calC (S, \xi )$. Here, we disregard subsurfaces of complexity $\leq \xi$ as per Remark \ref{heuristic}. The proof follows that of Theorems \ref{dfmc} and \ref{dfpc} above.

\begin{theorem}
    If $\mu,\nu\in \calC (S, \xi )$,   there exists a constant $K  (S ) $, 
    depending only on the topology of $S$, such that for each 
    $K>K  (S ) $ there exists $a\geq 1$ and $b\geq 0$ for which:
    $$
    d_\xi (\mu,\nu) \approx_{a,b} 
    \sum_{W\subseteq  S , \xi (W) > \xi } 
 Tf_K{d_{\calC (W)}}(\pi_{Y}(\mu),\pi_{Y}(\nu))
    $$
\label{dfic}
\end{theorem}

\section{Rank of $\calC (S, \xi )$}

\subsection{Lower Bound on Rank of $\calC (S, \xi )$}

In this section, we generalise a result of Brock-Farb \cite{brock-farb} to the case of interpolating graphs.
There are two extra ingredients. First, the observation in Remark \ref{ic} that $\calC (W, \xi )$ is quasi-isometric to the curve complex $\calC (W)$ if $\xi (W) = \xi +1$. The other ingredient is Theorem \ref{dfic} which generalises Theorem \ref{dfmc} and Theorem \ref{dfpc}.

Recall that a {\em quasiflat} in a metric space $X$ is a quasi-isometric embedding of Euclidean $n$-space in $X$; also, (Gromov \cite{gromov-ai} Section 6.2) that the {\em rank} of a metric
space $X$ is the maximal dimension $n$ of a {\em quasi-flat} in $X$.

As in \cite{brock-farb},
we shall say that  $S$ {\em decomposes} into essential subsurfaces
$R_1, \ldots, R_k$ if each $R_j$ is essential and if  $R_1, \ldots, R_k$ may be modified by an isotopy so that
they are pairwise disjoint and $S - R_1 \sqcup \ldots \sqcup
R_k$ is a collection of open annular neighborhoods of simple closed
curves on $S$, each isotopic to a boundary component of some $R_j$.

Fix $\xi \in \bN$. Let $r_\xi   (S ) $ denote the maximum number $k$ such that $S$ decomposes
into essential subsurfaces $R_1,
\ldots, R_k, T$ such that each $R_j$,
$j=1,\ldots, k$ has $\xi (R_j) = \xi + 1$ and $T$ is either empty or has $\xi (T) \leq \xi$.

\begin{theorem}
The graph $\calC (S, \xi )$ contains a quasi-flat of dimension $r_\xi   (S ) $.
\label{lb}
\end{theorem}

{\bf Proof:}  The proof is essentially a reworking 
in the present context of complexity $\xi$ graphs of Theorem 4.2 of \cite{brock-farb} by Brock and Farb.
By definition of $r_\xi = r$(say), the surface $S$ decomposes into subsurfaces
$$R_1, \ldots, R_{r(S)}, T$$ so that $\xi (R_j) = \xi + 1$ for each
$j$ and either $T$ is empty or $\xi (T) \leq \xi$.  

We now construct a quasi-isometric
embedding of the Cayley graph for $\bZ^r$ with the standard
generators into the complex $\calC (S, \xi )$.  

Let $\{ c_{j} \}$ be a   pants decomposition of $R_j$ and $\gamma$   a pants decomposition  of $T$. Along with the core curves of the open annuli
in $S - R_1 \sqcup \ldots \sqcup R_r \sqcup T$, 
the curves {\bf base}$c_{j}$ and {\bf base}$\gamma$
form a pants decomposition of $S$.

We let $g_j \colon \bZ \to \calC(R_j , \xi)$ be a (bi-infinite) geodesic so that
$g_j(0) = c_{j}$.  Since $\calC(R_j , \xi)$ is quasi-isometric to $\calC (R_j )$ (Remark \ref{ic} ) we might as well assume that $g_j \colon \bZ \to \calC(R_j  )$ is a quasigeodesic (though strictly speaking, we should compose $g_j$ with the quasi-isometry between $\calC(R_j , \xi)$ and  $\calC (R_j )$). Note that the quasi-sometry constants depend only on $R_j$ and hence only on the topology of $S$. Further, we identify $\calC(R_j , \xi)$ and  $\calC (R_j )$ via these uniform quasi-isometries. For the rest of this proof, we shall assume that $\calC(R_j , \xi) = \calC (R_j )$ rather than just being quasi-isometric to it. 

We define  the embedding
$ Q \colon \bZ^r \to \calC (S , \xi ) $
as: \\
$$Q(k_1, \ldots, k_r)  = (g_1(k_1), \ldots, g_r(k_r))$$.

Let $\vec{k} = (k_1,\ldots, k_n)$ and $\vec{l} = (l_1, \ldots, l_n)$.
Since elementary moves in the pants graph and hence in $\calC (S, \xi )$ along $g_j$ can be made independently in each
$R_j$, we have
$$d_{\xi}(Q(\vec{k}), Q(\vec{l})) \le \sum_{j=1}^n 
|l_j - k_j|=d_{\bZ^r}(\vec{k},\vec{l})$$
which shows that $Q$ is $1$-Lipschitz.

Given $R_j$, the subsurface
projection $\pi_{R_j}(Q(\vec{k}))$ to $R_j$  simply
picks out the curve $g_j(k_j)$ so we have
$$\pi_{R_j}(Q(\vec{k})) = g_j(k_j).$$ 
Thus, the projection distance
$$d_{R_j}(Q(\vec{k}), Q(\vec{l})) = 
d_{R_j}(g_j(k_j), g_j(l_j)) =  |k_j - l_j|$$ 
 since we are identifying $\calC (R_j , \xi )$ and 
$\calC(R_j)$.

By Theorem \ref{dfic}, there exists $K_0=K_0  (S ) $ so that for all
$K\geq K_0$ there exist constants $a$ and $b$ 
so that if we let $P_{\vec{k}} = Q(\vec{k})$ and 
$P_{\vec{l}} =
Q(\vec{l})$ then we have the inequality $$ \sum_{\stackrel{Y \subseteq
S \, \xi(Y) > \xi}{d_{\calC (Y)}(\pi_Y(P_{\vec{k}}), \pi_Y(P_{\vec{l}})) >M}} 
d_Y(P_{\vec{k}},P_{\vec{l}})
\le a \, d_{\xi}(P_{\vec{k}}, P_{\vec{l}}) + b.$$
But the left-hand-side of the inequality is bounded below by
$$\max_j |k_j - l_j| \ge \frac{\sum_j |k_j - l_j|}{r}.$$
Thus, $Q$ is a quasi-isometric embedding.
$\Box$

\subsection{Upper Bound on Rank of $\calC (S, \xi )$}

In this section, we generalise a recent result of Behrstock-Minsky \cite{behrstock-minsky} to the context of interpolating graphs. As we shall be following 
\cite{behrstock-minsky} closely, we shall indicate only the steps in the argument and the necessary modifications. 

\medskip

\noindent
 {\bf \underline{Step 1: Coarse Product Regions in $\calC (S, \xi )$}}

\smallskip 

\noindent
First, we describe the geometry of the set of pants decompositions in $\calC (S, \xi )$
containing a prescribed set of base curves. Equivalently, in the coned-off 
mapping class group $\Gamma (S, \xi )$, such a set corresponds to the coned off coset of the
stabilizer of a simplex in the complex of curves. 
These regions coarsely decompose as products as in \cite{behrstock-minsky}.

Let $\Delta$ be a  multicurve
 in $S$.  Partition $S$ into
subsurfaces isotopic to complementary components of $\Delta$. Throw away
components homeomorphic to $S_{0,3}$. Obtain the
 {\it partition} of $\Delta$,  denote it as $\sigma(\Delta)$.

Theorem \ref{dfic}  gives the
following generalisation of Lemma 2.1 of Behrstock-Minsky \cite{behrstock-minsky}. 

\begin{lemma}
Let $\calQ (\Delta)\subset \calC (S, \xi )$ denote the set of pants decompositions whose bases
contain $\Delta$. Then, sending each (family of base curves of a) pants decomposition to the restrictions to elements of $\sigma (\Delta )$ we obtain a quasi-isometric identification 
\begin{equation*}\label{Q product}
\calQ(\Delta) \approx \prod_{U\in \sigma(\Delta), \xi (U) > \xi} \calC (U, \xi )
\end{equation*}
with uniform constants.
\label{qipdkt}
\end{lemma}

\medskip

\noindent
 {\bf \underline{Step 2: Ultralimits of $\calQ (\Delta)$}}

\smallskip 

\noindent
We refer the reader to \cite{behrstock-minsky} Section 1.4 for the necessary background on ultralimits and asymptotic cones.
$ \AC(S, \xi)$ will denote an asymptotic cone of $\calC (S, \xi )$ and $\mu_0$ a preferred base-point. 

\begin{defn}
For a  sequence $\seq\Delta=\{\Delta_n\}$ 
such that $\ulim \frac{1}{s_n}d_\xi (\mu_0, \calQ(\Delta_n)) < \infty$, define $\AQ (\seq\Delta) \subset \AC (S, \xi )$ to be the
ultralimit of $\calQ (\Delta_n)$, with  metrics rescaled by $1/s_n$.
\end{defn}

Also, since the topological type of $\sigma(\Delta_n)$ is $\omega$-a.e.\ constant, we may define $\sigma(\seq\Delta)$ to be the $\omega$-limit of $\sigma(\Delta_n)$'s. Note also that the complexities $\xi (U_{nj})$ are 
$\omega$-a.e.\ constant for components $U_{nj}$ of $U_n$.

Then Lemma \ref{qipdkt} and the fact that ultralimits commute 
with finite products gives the following generalisation of Equation 2.2 of \cite{behrstock-minsky}:

\begin{lemma} There is a uniform
bi-Lipschitz identification 
\begin{equation*}
\AQ(\seq\Delta) \homeo 
\prod_{\seq U \in \sigma(\seq\Delta) , \xi (\seq U ) > \xi} \AC(\seq
U, \xi). 
\end{equation*}
\label{blpdkt}
\end{lemma}

\medskip

\noindent
 {\bf \underline{Step 3: $\bR$-trees and Product Regions in Asymptotic Cones}}

\smallskip 

\noindent
The following definition is adapted from  Behrstock \cite{behrstock-asymptotic} and Behrstock-Minsky \cite{behrstock-minsky}.

\begin{defn} 
Let $\seq W=(W_n)$ be a sequence of connected subsurfaces (considered
mod $\omega$) and $\seq x\in \AC(\seq W , \xi )$,  
Then  $F_{\seq W,\seq x , \xi }\subset \AC(\seq W , \xi )$ is defined as:
$$
F_{\seq W,\seq x , \xi } = \{ \seq y\in\AC(\seq W , \xi ) : d_{\AC(\seq W , \xi )}(\seq x,\seq y)=0\ 
\  \text{for all  proper subsets $\seq U\subset \seq W$} \}.
$$
\end{defn}

The next theorem is a version of a theorem due to  Behrstock \cite{behrstock-asymptotic} adapted to our context of interpolation graphs.

\begin{theorem}
Let $\seq W = (W_n)$ be a sequence of connected proper subsurfaces of 
$S$, and
$\seq x \in \AC (\seq W, \xi )$. 
Any two points $\seq y,\seq z\in F_{\seq W,\seq x, \xi}$ are connected by a
unique embedded path in 
$\AC(\seq W , \xi )$, and this path lies in 
$F_{\seq W,\seq x}$. 
In particular,  $F_{\seq{W,x} , \xi}$ is an $\bR$-tree. 
\label{r-tree}
\end{theorem}

Next, for  $\seq W $ and $\seq x$ as above, 
separating product regions in $\AC(\seq W , \xi )$, denoted  $P_{\seq W,\seq x , \xi}$, are
subsets of $\AQ(\partial \seq W)$ defined as follows: \\
In the bi-Lipschitz product structure on $\AQ(\partial \seq W)$ (Lemma \ref{blpdkt}), 
 $\seq W$ is a member of $\sigma(\partial 
\seq W)$. Therefore, $\AC (\seq W , \xi )$ 
appears as a factor. 
Define $P_{\seq W,\seq x, \xi}$
 to be the subset of
$\AQ(\partial \seq W , \xi )$ consisting of points whose coordinate in
the $\AC(\seq W, \xi)$ factor lies in $F_{\seq W,\seq x , \xi}$. The following Lemma generalises Lemma 3.3 of \cite{behrstock-minsky}.

\begin{lemma}
There exists a  bi-Lipschitz identification of
$P_{\seq W,\seq x , \xi}$ with
$$
F_{\seq W,\seq x , \xi} \times \AC(\seq W^c , \xi).
$$
\label{pset}
\end{lemma}

\medskip

\noindent
 {\bf \underline{Step 4: Global Projection Maps}}

\smallskip 

\noindent
The following Theorem generalises Theorem 3.5 of \cite{behrstock-minsky} and gives a global projection map for $F_{\seq W, \seq x , \xi}$: 

\begin{theorem}
Given $\seq x\in\AC(\seq W , \xi )$, there is a continuous map
$$
\Phi=\Phi_{\seq W, \seq x , \xi} \co \AC(S, \xi ) \to F_{\seq W, \seq x , \xi }
$$
with these properties:  
\begin{enumerate}
\item 
$\Phi$ restricted to $P_{\seq W, \seq x , \xi}$ 
is projection to the first factor in the product structure 
$P_{\seq W, \seq x , \xi} \homeo
  F_{\seq W, \seq x , xi}\times \AC(\seq W^c , \xi )$.
\item
$\Phi$ is locally constant in the complement of $P_{\seq W, \seq x , xi}$. 
\end{enumerate}
\label{globalproj}
\end{theorem}

\medskip

\noindent
 {\bf \underline{Step 5: Separating Sets}}

\smallskip 

\noindent The sets $P_{\seq W, \seq x , \xi} \homeo
  F_{\seq W, \seq x , xi}\times \AC(\seq W^c , \xi )$ give rise to separating sets in $\AC(S , \xi )$ as in Theorem 3.6 of \cite{behrstock-minsky}.

\begin{theorem}
There is a family $\calL$ of closed subsets of $\AC (S, \xi )$ such that
any two points in $\AC (S, \xi )$ are separated by some $L\in\calL$. 
Moreover each $L\in \calL$ is isometric to
$\AC (Z, \xi )$, where $Z$ is some proper essential (not necessarily
connected) subsurface of $S$, with $r(Z) < r(S)$. 
\label{separate}
\end{theorem}

\medskip

\noindent
 {\bf \underline{Step 6: Inductive Dimension}}

\smallskip 

\noindent 
Using the separating sets, we complete the argument as in Theorem 4.1 of \cite{behrstock-minsky}. $\ind$, $\Ind$ and $dim$ will denote small inductive dimension, $\ind$,
large inductive dimension,  and covering (or topological) dimension respectively. 

$\indhat(X)$ denote the supremum of $\ind(X')$ over all
locally-compact subsets $X'\subset X$;  similarly  $\Indhat$
and $\dimhat$. Then using Theorem \ref{separate} above the argument by Behrstock and Minsky yields

\begin{theorem}
$\indhat(\AC(S, \xi ))=\Indhat(\AC(S, \xi ))=\dimhat(\AC(S, \xi )) =  r_\xi (S)$.
\label{ind}
\end{theorem}

Since a quasiflat of dimension $r$ would yield a dimension $r$ locally compact subset of $\AC(S, \xi )$, Theorem \ref{ind} immediately gives:

\begin{cor}
If the graph $\calC (S, \xi )$ contains a quasi-flat of dimension $r$, then $r \leq r_\xi   (S ) $.
\label{ub}
\end{cor}

Combining Theorem \ref{lb} and Corollary \ref{ub}
we obtain the main theorem of this paper:

\begin{theorem}
The rank of the interpolating graph $\calC (S, \xi )$ of complexity $\xi$ is $r_\xi   (S ) $.
\label{rank}
\end{theorem}

\subsection{Problems}
A number of (hopefully interesting) issues arise from the notion of interpolating graphs of complexity $\xi$.

\medskip

\noindent {\bf Conjecture1: Rank one implies Hyperbolic}\\
If for $\xi > 0$,  $r_\xi (S) = 1$  
 is $\calC (S, \xi )$ a hyperbolic metric space?

\smallskip

This  is motivated by \\

\begin{enumerate}
\item the observation (Remark \ref{ic} ) that $\calC (S, \xi (S) -1 )$ is quasi-isometric to the curve graph and that in this case, $r_{\xi (S) -1}$ is clearly one (all of 
$S$ is embedded in $S$ by the identity map).\\
\item
Next, for $\xi = -2$, $r_\xi = 1$ means that the maximum number of disjoint homotopically distinct annulii (connected subsurface of complexity $(-2 + 1 =1)$) that $S$ admits is precisely one. Then $\xi (S) =1$ and it is precisely these cases that have hyperbolic marking complex $\calC (S, -2)$. \\
\item
$\xi = -1$ corresponds to the pants graph. This case is special and here the appropriate hypothesis would be  $r_0(S) = 1$ because $S_{0,3}$ has trivial curve complex; hence in order to get an infinite diameter curve complex, we have to step up $\xi = -1$ by $2$ and then calculate the maximum number of disjoint embedded subsurfaces of complexity $\geq 1$.
Brock and Farb \cite{brock-farb} show that the pants graph is hyperbolic iff $\xi (S) = 2$ and each of the two possibilities (5-holed sphere and two-holed torus) admit exactly one (disjoint) subsurface of complexity $1$.
\item
$\xi = 1$, $\xi (S) = 3$. Recent work of Brock and Masur \cite{brock-masur-cwp} is closely related to this case.
\end{enumerate}

Thus the above Conjecture would serve to unify all the above cases.

\medskip

\noindent {\bf Problem 2}
Is the automorphism group of the interpolating graph commensurable
with the mapping class group $MCG(S)$ ?

\smallskip

This question is a special case of 
Ivanov's metaconjecture \cite{ivanov-15} that every object naturally associated to a surface $S$ and hav-
ing a suﬃciently rich structure has $MCG(S)$ as its groups of automorphisms.

One other piece of motivation is Margalit's result \cite{margalit-pants} reducing the automorphism group of the pants graph to that of the curve complex  C(S).

\medskip

\noindent {\bf Problem 3} There is a hierarchy
$$\calC (S, -2) \rightarrow \calC (S, -1) \rightarrow \calC (S, 1) \cdots C(S, \xi ) \rightarrow C(S, \xi + 1) \rightarrow \cdots \rightarrow C(S, \xi (S) - 1)$$
where the first term is the marking graph, the second the pants
graph, the last the curve graph. 

The map $C(S, \xi ) \rightarrow C(S, \xi + 1)$ is given by coning
a collection of subsets corresponding to curve graphs of subsurfaces of complexity
$\xi + 1$. Thus the last term is hyperbolic and the preimages
of points at each stage are hyperbolic by Masur-Minsky's Theorem \ref{cc-hyp}. 
This raises the hope that the hierarchy paths constructed by Masur and Minsky in \cite{masur-minsky2} may alternately be  inductively constructed in a bottom-up approach from the curve complex. Further, at each stage we should obtain a hierarchy path in the interpolating graph $\calC(S, \xi )$. 

\medskip

\noindent {\bf Problem 4} (independently due to Yair Minsky \cite{minsky-pc} ) Finally, there ought to be a general geometric structure lying between strong and weak relative hyperbolicity of which the mapping class group is a special case. Let us call this putative structure {\em graded relative hyperbolicity}
(terminology independently due to Yair Minsky \cite{minsky-pc} ) . A possible definition would be the existence of a sequence 
$$X_n \rightarrow X_{n-1} \rightarrow  \cdots X_1$$
of spaces and maps where at the $i$th stage one cones off a collection $\calC_i$ of (uniformly) quasiconvex hyperbolic subsets of $X_i$. Further, we demand that $X_1$ be hyperbolic. 

A toy example is given by $X = X_n = \Gamma_G$ the Cayley graph of a hyperbolic group $G$. $H$ is assumed to be a quasiconvex subgroup of height $n$ (see Gitik-Mitra-Rips-Sageev \cite{GMRS} for instance). In passing from $X_i$ to $X_{i-1}$, we cone all cosets of $\cap_{j = 1 \cdots i-1} g_j H g_j^{-1}$ for essentially distinct cosets $g_jH$. In this particular case, the role of hierarchy paths might be taken by {\em electro-ambient quasigeodesics} introduced by the author in \cite{brahma-ibdd}. 

\bibliography{intcx}
\bibliographystyle{plain}

\end{document}